\documentclass{article}

\usepackage{amsmath}

\title{%
      Connection formulas
  for the confluent hypergeometric functions
  and the functional relation
  for the Hurwitz zeta function}
\author{Michitomo NISHIZAWA  and Kimio UENO}
\date{}

\begin{document}
\maketitle
\begin{abstract}
The functional relation of the Hurwitz zeta function is proved by using 
the connection problem of the confluent hypergeometric equation. 
\end{abstract}

\section{Introduction}

In this article, we give an alternative proof for the functional relation of
the Hurwitz zeta function by using the connection problem for
the confluent hypergeometric equation. 
This proof appeared first in \cite{UN1} in an implicit form, and second in \cite{UN2} 
which was written unfortunately in Japanese. In these articles we considered mainly
the functional relation for the q-analogue of the Hurwitz zeta function.
So we decided to publish again the proof in a complete form. \\

\noindent
\textbf{Acknowledgement} \quad The authors would like to express their deep gratitude
to the organizing committee of the conference ZTQ, especially to 
Professor Takashi Aoki and Professor Yasuo Ohno
for giving the authors the opportunity presenting their research.  

The second author is
partially supported by JPSP Grant-in-Aid No. 15540050 and 
Waseda University Grant for Special Research Project (2002A-067,
2003A-069).

\section{Hurwitz zeta function}

The generalized zeta function introduced by Hurwitz (Hurwitz zeta function, for short)
is, by definition, 
\begin{eqnarray}
   \zeta(s,z) \,=\, \sum_{k=0}^{\infty}\frac{1}{(k+z)^s},
\end{eqnarray}
where we suppose that $0 \leq z < 1$. The series is absolutely convergent for $\Re s > 1$,
and is analytically continued to the whole $s$-plane as a meromorphic function. 
Evidently, $\zeta(s,1) \,=\, \zeta(s)$, which is the Riemann Zeta function. 
Furthermore, they satisfy a functional relation
\begin{eqnarray}\label{HURWITZ}
    \zeta(s,z) \,=\, \Gamma(1-s) \left\{(2\pi i)^{s-1}\mathcal{L}(1-s,z)
                      +(-2\pi i)^{s-1}\mathcal{L}(1-s,1-z) \right\},
\end{eqnarray}
which was established by Hurwitz himself \cite{WW}. Here 
\begin{eqnarray*}
      \mathcal{L}(s,z) \,=\, \sum_{n=1}^{\infty}\frac{e^{2\pi inz}}{n^s}
\end{eqnarray*}
is a generalized polylogarithm and $\Gamma(s)$ is Euler's gamma function.
Let call (\ref{HURWITZ}) the Hurwitz relation.
When $z=1$, it reduces to the functional equation for the Riemann zeta function:
\begin{eqnarray*}
  \zeta(s) \,=\, 2^s\pi^{s-1}\Gamma(1-s)\sin\left(\frac{\pi s}{2}\right)\zeta(1-s).
\end{eqnarray*}
To investigate analytic continuation of zeta functions, the sum formula of Euler-Maclaurin
\begin{eqnarray}\label{EM}
     \sum_{r=0}^{\infty}f(r) \,=\, \int_0^{\infty}\!\!f(t)\,dt
       +\sum_{k=1}^n \, \frac{B_k}{k!}\{ f^{(k-1)}(\infty)-f^{(k-1)}(0) \} \nonumber \\
       +\frac{(-1)^n}{n!}\int_0^{\infty}\!\!\overline{B}_n(t)f^{(n)}(t)\,dt
\end{eqnarray}
is a useful tool. Here $\overline{B}_n(t)=B_n(t-[t])$ ($[t]$ denotes the integral 
part of $t$),
$B_n(t)$ is the n-th Bernoulli polynomial defined by
\begin{eqnarray*}
         \sum_{n=0}^{\infty}\frac{B_n(t)u^n}{n!} \,=\,
			\frac{{}ue^{tu}{}}{e^u-1},
\end{eqnarray*}
and $B_n=B_n(0)$ is the n-th Bernoulli number. Putting $f(t)=(t+z)^{-s}$, 
and $n=2$ in (\ref{EM}), we have
\begin{eqnarray}\label{EM1}
  \zeta(s,z) \,=\, \frac{z^{1-s}}{s-1}+\frac{z^{-s}}{2}+\frac{sz^{-s-1}}{12}
            - \int_0^{\infty}\frac{\overline{B}_2(t)}{2!}
                \left(\frac{d}{dt}\right)^2\left\{\frac{1}{(z+t)^s}\right\}\,dt
\end{eqnarray}
Substituting the Fourier expansion
\begin{eqnarray}
       \overline{B}_2(t)=\frac{1}{\pi^2}\sum_{n=0}^{\infty}\frac{\cos(2\pi nt)}{n^2}
\end{eqnarray}
into (\ref{EM1}), making once partial integration, we obtain
\begin{eqnarray}\label{AN}
  \zeta(s,z) \,=\, \frac{z^{1-s}}{s-1}+\frac{z^{-s}}{2}+\frac{sz^{-s}}{2\pi i}
        \sum_{l\neq0}\frac{1}{\,l\,}
                 \int_0^{\infty}\frac{e^{-2\pi il zu}}{(1+u)^{s+1}}\,du.
\end{eqnarray}

\section{Confluent hypergeometric function}

Now we observe the confluent hypergeometric equation
\begin{eqnarray}\label{CONFL}
   x\frac{d^2y}{dx^2}+(\gamma-x)\frac{dy}{dx}-\alpha y = 0.
\end{eqnarray}
This equation has a regular singular point at $x=0$, and an irregular singular point
at $x=\infty$. The confluent hypergeometric series
\begin{eqnarray}\label{CONFLFUNC}
    F(\alpha,\gamma;x) \,=\, \sum_{n=0}^{\infty} 
             \frac{(\alpha)_nx^n}{(\gamma)_nn!},
\end{eqnarray}
and \ $x^{1-\gamma}F(\alpha-\gamma+1,2-\gamma;x) \, (\,=e^xx^{1-\gamma}F(1-\alpha,2-\gamma;-x))$
\ form a system of fundamental solutions to (\ref{CONFL}) around $x=0$.
Let us introduce a function defined by
\begin{eqnarray}\label{U}
  U(\alpha,\gamma;x) \,=\, \frac{1}{\Gamma(\alpha)}
                     \int_0^{\infty}e^{-xu}(1+u)^{\gamma-\alpha-1}u^{\alpha-1}\,du.
\end{eqnarray}
This is a solution to (\ref{CONFL}) around $x=\infty$ \cite{S}, 
and is connected to the former solutions by
\begin{eqnarray}\label{CONN}
     U(\alpha,\gamma;x) \,=\, \frac{\Gamma(1-\gamma)}{\Gamma(1+\alpha)}F(\alpha,\gamma;x)
         + \frac{\Gamma(\gamma-1)}{\Gamma(\alpha)}e^xx^{1-\gamma}F(1-\alpha,2-\gamma;-x)).
\end{eqnarray}
Furthermore, for $\Re \alpha>0$, we see that
\begin{eqnarray}\label{ASYM}
        U(\alpha,\gamma;x) \,\sim\, x^{-\alpha}, \qquad (|x| \to \infty).
\end{eqnarray}
Substituting (\ref{U}) to (\ref{AN}) with $\alpha=1, \gamma=1-s, \text{and} \ x=-2\pi il z$,
we obtain
\begin{eqnarray}\label{EM2}
  \zeta(s,z) \,=\, \frac{z^{1-s}}{s-1}+\frac{z^{-s}}{2}+\frac{sz^{-s}}{2\pi i}
        \sum_{l\neq0}\frac{1}{\,l\,}U(1,1-s;-2\pi i l z).
\end{eqnarray}
Here we should note that the infinite sum in (\ref{EM2}) is absolutely convergent due to 
the asymptotic behavior (\ref{ASYM}). According to (\ref{CONN}), 
we have a connection formula
\begin{eqnarray}
U(1,1-s;-2\pi il z) \,=\, \frac{1}{s}\left\{ F(1,1-s;-2\pi il z) - 
                              \Gamma(1-s)(-2\pi il z)^se^{-2\pi il z} \right\}.
\end{eqnarray}
Hence
\begin{eqnarray}
  \frac{sz^{-s}}{2\pi i}\sum_{l\neq0}\frac{1}{\,l\,}U(1,1-s;-2\pi i l z) = 
  \frac{z^{-s}}{2\pi i}\sum_{l\neq0}\frac{1}{\,l\,}F(1,1-s;-2\pi i l z)+ \nonumber \\
  + \Gamma(1-s) \left\{(2\pi i)^{s-1}\mathcal{L}(1-s,z)
                      +(-2\pi i)^{s-1}\mathcal{L}(1-s,1-z) \right\}.
\end{eqnarray}
By virtue of the integral representation for $F(\alpha,\gamma;x)$
\begin{eqnarray*}
  F(\alpha,\gamma;x) \,=\, \frac{\Gamma(\gamma)}{\Gamma(\gamma-\alpha)\Gamma(\alpha)}
        \int_0^1 e^{zt}t^{\alpha-1}(1-t)^{\gamma-\alpha-1}\,dt,
\end{eqnarray*}
we have
\begin{eqnarray*}
\frac{1}{2\pi i}\sum_{l\neq0}\frac{1}{\,l\,}F(1,1-s;-2\pi i l z) &=&
\frac{s}{\pi} \sum_{n=1}^{\infty} \int_0^1 \frac{\sin(2\pi nzt)}{n}(1-t)^{-s} \,dt \nonumber \\
&=& \frac{s}{\pi^2z} \left\{\frac{\pi^2}{6}+(s+1)
                \int_0^1 \overline{B}_2(t)(1-t)^{-s-2}\,dt \right\}.
\end{eqnarray*}
Since $\overline{B}_2(t)=\pi^2(t^2-t+1/6) \quad (0 \leq t \leq 1)$, we have
\begin{eqnarray*}
\frac{z^{1-s}}{s-1}+\frac{z^{-s}}{2}+\frac{z^{-s}}{2\pi i}
             \sum_{l\neq0}\frac{1}{\,l\,}F(1,1-s;-2\pi i l z) \,=\, 0
\end{eqnarray*}
for $0 \leq z <1$. Thus we obtain the Hurwitz relation (\ref{HURWITZ}).

In the case of the q-analogue of the Hurwitz zeta function, the connection formulas for the
hypergeometric function play the same role as in the present case \cite{UN1}.

\vspace{5mm}

\begin{tabular}{l}
 NISHIZAWA Michitomo\\
 Graduate School of Mathematical Sciences\\
 University of Tokyo\\
 Komaba 3-8-1, Meguro-ku, Tokyo 153-8914, Japan\\
 E-mail: {\ttfamily mnishi@ms.u-tokyo.ac.jp}
\end{tabular}

\begin{tabular}{l}
 UENO Kimio\\
 Department of Mathematical Sciences\\
 School of Science and Engineering\\
 Waseda University\\
 Okubo 3-4-1, Shinjuku-ku, Tokyo 169-8555, Japan\\
 E-mail: {\ttfamily uenoki@mse.waseda.ac.jp}
\end{tabular} 

\end{document}